\documentstyle[12pt,amssymb,amstex]{amsart}
\numberwithin{equation}{section}

\def\ch{{\cal H}}

\def\ck{{\cal K}}

\def\co{{\cal O}}

\def\cq{{\cal Q}}

\def\cx{{\cal X}}

\def\ga{{\frak A}}
\def\gb{{\frak B}}

\def\gd{{\frak D}}

\def\gf{{\frak F}}

\def\gi{{\frak I}}

\def\gk{{\frak K}}

\def\gam{{\frak M}}

\def\bc{{\Bbb C}}

\def\bk{{\Bbb K}}
\def\bm{{\Bbb M}}
\def\bn{{\Bbb N}}

\def\a{\alpha}
\def\b{\beta}
  
  \def\D{\Delta}
\def\eps{\varepsilon}

\def\r{\rho}
\def\s{\sigma}

\def\f{\varphi} \def\F{\Phi}
\def\th{\theta}
\def\om{\omega} \def\Om{\Omega}

\def\itm#1{\item[($#1$)]}

\def\id{\hbox{id}}
\def\ker{\hbox{Ker}}

\newtheorem{thm}{Theorem}
\newtheorem{lem}{Lemma}

\newtheorem{prop}{Proposition}

 \theoremstyle{definition}
\newtheorem{defin}{Definition}

 \theoremstyle{remark}
\newtheorem{rem}{Remark}

\begin{document}

\title[operator ideals]
{Some operator ideals in non-commutative functional analysis} 
\author{Francesco Fidaleo}
\address{Francesco Fidaleo\\
Dipartimento di Matematica\\
II Universit\`{a} di Roma (Tor Vergata)\\
Via della Ricerca Scientifica, 00133 Roma, Italy}
\email{{\tt Fidaleo@@mat.utovrm.it}}
\maketitle
 \begin{abstract}
We characterize classes of linear maps between operator spaces $E$, $F$
which factorize through maps arising in a natural manner via 
the Pisier vector-valued non-commutative $L^p$
spaces $S_p[E^*]$ based on the Schatten classes on the separable Hilbert 
space $l^2$. These classes of maps can be viewed as quasi-normed operator
ideals in the category of operator spaces, that is in non-commutative
(quantized) functional analysis. The case $p=2$ provides a Banach
operator ideal and allows us to characterize the split property
for inclusions of $W^*$-algebras by the $2$-factorable maps. The various 
characterizations of the split property have interesting applications
in Quantum Field Theory.
 \vskip 0.3cm \noindent
Mathematics Subject Classification Numbers: Primary 47D15, 47D25, 
Secondary 46L35.\\
Key words: Linear spaces of operators,
Operator algebras and ideals on Hilbert spaces; Classifications, Factors.
 \end{abstract}

\section{Introduction}
 
The study of classes of operator ideals in Hilbert and Banach space
theories has a long history. First some interesting classes of
maps between classical function spaces were considered 
and therefore many classes of operators were so intensively studied. 
Due to the wideness
of the subject we refer the reader to the monoghraphies 
\cite{G1,G2,J,Pi,Pi1} and
to the references quoted therein. Many of these classes of maps
can be considered as operator ideals in the category of Banach
spaces as is well exposed in \cite{Pi1} where a wide class of operator
ideals is treated from an axiomatic as well as concrete viewpoint. 
On the other
hand, if a Banach space $E$ is equipped with a sequence
of compatible norms on all the matrices $\bm_n(E)$ with entries in $E$, 
it can be viewed as a subspace of a $C^*$-algebra that is a
(concrete) operator space \cite{R}. In this context the natural
arrows between operator spaces are the completely bounded linear
maps. These ideas can be interpreted as a non-commutative
(quantized) version the functional analysis as is well explained in
several papers, see e.g. \cite{B,E,E1,ER6,Pa1,R,We}. 
Also in the operator spaces context,
interesting classes of completely bounded linear maps have been introduced and
studied, see \cite{ER3,ER4,ER5,F,1P,P2,P4}. 
These classes of maps can be naturally considered
as operator ideals between operator spaces i.e. in non-commutative
functional analysis. Namely one should replace 
the bounded maps, i.e. the natural arrows between Banach spaces, 
in the definition of an operator ideal given in \cite{Pi1}, 
with the completely bounded maps which are the natural arrows 
in the category of operator spaces.\\
In this paper we characterize, for each $1\leq p<+\infty$, 
classes $\gf_p(E,F)$ of linear maps between operator 
spaces $E$, $F$ which factorize through maps arising in a natural manner via 
the Pisier vector-valued non-commutative $L^p$
spaces $S_p[E^*]$ based on the Schatten classes on the separable Hilbert 
space $l^2$. We call these maps the (operator) $p$-{\it factorable} maps and
omit here the word ``operator" as in the following there is no matter
of confusion with the factorable operators of Banach space context 
\cite{Kw,Pi1} (see the end of Section 6 below for some comments about 
this point).\\
The classes of $p$-factorable maps $\gf_p$ are quasi-normed (complete) 
operator ideals in the operator space setting.
Moreover, there is a geometrical description for
the image of the unit ball $T(E_1)\subset F$
under a $p$-factorable injective map, then, under certain conditions on
the image space $F$. Such
a description of the shape of $T(E_1)$ is analogous to those
arising from the Banach space context, see \cite{B1};
or from the operator space context, see \cite{F}.
The case with $p=2$ is particular: $\gf_2$ is a Banach operator
ideal; furthermore one can give the geometrical description for 
the image of the unit ball $T(E_1)\subset F$, $T\in\gf_2$ injective,
without adding any other condition on the image space $F$. 
Finally, as an application,
we outline a characterization of the split property for certain inclusion
of $W^*$-algebras (\cite{B1,D,DL2,F}) and,
in particular, for inclusions of von Neumann algebras of local observables
arising from Quantum Field Theory as described in\cite{B2,B3,B5,Su}.\\

This paper is organized as follows.\\

After some preliminaries about the operator 
spaces, we discuss some interesting result relative to 
non-commutative vector-valued $L^p$ spaces recently introduced 
by Pisier in \cite{P2,P4} which are of main interest here.
Successively we define, for each $1\leq p<+\infty$, a class of linear maps
$\gf_p(E,F)$ between the operator spaces $E$, $F$ which are limits
of finite rank operators,  hence compact ones. We call such classes of maps
the ideals of $p$-{\it factorable} maps. The spaces
$\gf_p(E,F)\subset \gk(E,F)$ are in a natural way operator ideals between
operator spaces. If $p=2$ we obtain a Banach operator ideal.
A further section is devoted to provide a geometrical description of 
$T(E_1)\subset F$ where $T\in\gf_p(E,F)$ is an injective operator. 
The description of the shape of $T(E_1)$ 
can be made for a general $p$ under certain conditions on the image
space $F$ (i.e. if $F$ is an injective $C^*$-algebra) or 
without adding any other condition on $F$ if $p=2$.\\ 
As an application, we conclude with a section where we outline a 
characterization of the split property for inclusion $N\subset M$ 
of $W^*$-factors 
in terms of $2$-factorable maps which is also applicable to inclusions of 
von Neumann algebras (which have a-priori a non-trivial center)  
of local observables arising in Quantum Field Theory. 
The detailed exposition of the complete characterization of the split 
property in terms of properties of the canonical non-commutative $L^2$
embedding ${\displaystyle N\stackrel{\F_2}{\hookrightarrow}L^2(M)}$ will be 
given in a separated paper \cite{F1}.

\section{On operator spaces}

For the reader's convenience we collect 
some preliminary results about the operator
spacesof which we need in the following. Details 
and proofs can be found in the sequel.\\
In this paper all the operator spaces are complete as normed spaces if it is
not otherwise specified.

\subsection{Operator spaces}

For an arbitrary normed space $X$, $X_1$ denotes its (closed) unit ball. We
consider a normed space $E$ together with a sequence of norms $\|\ \|_n$
on $\bm_n(E)$, the space of $n\times n$ matrices with entries in $E$. For
$a,b\in \bm_n$ these norms satisfy
\begin{align}
\begin{split}
&\|avb\|_n\leq\|a\|\|v\|_n\|b\|,\\
&\|v_1\oplus v_2\|_{n+m}=\max\{\|v_1\|_n,\|v_2\|_m\}
\label{ER}
\end{split}
\end{align}
where the above products are the usual row-column ones. This space with the above norms is called an
(abstract) operator space.\\

\vskip 0.1cm

If\ $T:E\to F$, $T_n:\bm_n(E)\to\bm_n(F)$ are given by $T_n:=T\otimes\id$. 
$T$ is said to be {\it completely bounded} if 
$\sup\|T_n\|:=\|T\|_{cb}<+\infty$; $\gam(E,F)$ denotes the set of all the 
completely bounded maps
between $E$, $F$. Complete contractions, complete isometries and complete
quotient maps have an obvious
meaning. It is an important fact (see \cite{R}) that a linear space $E$ 
with norms on each $\bm_n(E)$ has a realization as a concrete operator space 
i.e. a subspace of a $C^*$-algebra, if and only if these norms satisfy 
the properties in (\ref{ER}). We note that, if $dim(E)>1$ then there would
be a lot of non-isomorphic operator space structures on $E$, 
see e.g. \cite{Pa1}.\\

Given an operator space $E$ and $f=\bm_n(E^*)$, the norms
\begin{align}
\begin{split}
\|f\|_n:=\sup&\{\|(f(v))_{(i,k)(j,l)}\|:v\in\bm_m(E)_1,\ m\in\bn\}\\
(&f(v))_{(i,k)(j,l)}:=f_{ij}(v_{kl})\in \bm_{mn}
\label{dual}
\end{split}
\end{align}
determines an 
operator space structure on $E^*$ that becomes itself an operator space which
has been called in \cite{B} the {\it standard dual} of $E$.\\

We now consider the linear space $\bm_I(E)$ for any index $I$ as the $I\times I$
matrices with entries in $E$ such that
$$
\|v\|_{\bm_I(E)}:=\mathop{\sup}\limits_{\D}\|v^{\D}\|<+\infty
$$ 
where $\D$ denotes any finite 
subset of $I$. For each index set $I$,
$\bm_I(E)$ is in a natural way an operator space via the inclusion
$\bm_I(E)\subset\gb(\ch\otimes \ell^2(I))$ if E is realized as a subspace of
$\gb(\ch)$. Of interest is also the definition of $\bk_I(E)$ as those elements
$v\in\bm_I(E)$ such that $v=\lim\limits_{\D} v^{\D}$. Obviously
$\bm_I(\bc)\equiv\bm_I=\gb(\ell^2(I))$ and 
$\bk_I(\bc)\equiv\bk_I=\gk(\ell^2(I))$,
the set of all the compact operators on $\ell^2(I)$. For $E$ complete we remark
the bimodule property of $\bm_I(E)$ over $\bk_I$ because, for $\a\in\bk_I$,
$\a^{\D} v$, $v\a^{\D}$ are Cauchy nets in $\bm_I(E)$ whose limits
define unique elements $\a v$, $v\a$ that can be 
calculated via the usual row-column product.\\
Given an index set $I$ we can define as usual a map 
$\cx:\bm(E^*)\to\gam(E,\bm_I)$
which is a complete isometry, given by 
\begin{equation}
\label{mappa}
(\cx(f)(v))_{ij}:=f_{ij}(v).
\end{equation}
Moreover, if $f\in\bk(E^*)$, $\cx(f)$ is norm limit of finite rank
maps, then $\cx(\bk_I(E^*))\subset\gk(E,\bk_I)$ so $\cx(f)(v)\in\bk_I$
for each $v\in E$, see \cite{ER2}, pag 172.\\

\vskip 0.1cm

Remarkable operator space structures on a Hilbert space $H$ have been
introduced and studied in \cite{ER3} where natural identifications between
Banach spaces have been considered.
The following identifications 
\begin{align*}
&\bm_{p,q}(H_c):=\gb(\bc^q,H^p),\\
&\bm_{p,q}(H_r):=\gb(\overline{H}^q,\bc^p)
\end{align*}
define on $H$ the column and row structures respectively.\\

\vskip 0.1cm

Recently the theory of complex interpolation (\cite{BL}) has been developped
by Pisier in the context of operator spaces as well. An operator space
structure on Hilbert spaces has been so introduced and studied, that is 
the Pisier $OH$ structure, see \cite{P1,P3}. It can be viewed 
as an interpolating structure between $H_c$, $H_r$:
$$
OH(I):=(H_c,H_r)_{1/2}
$$
where the cardinality of the index set $I$ is equal to the (Hilbert)
dimension of $H$. The characterization of the $OH$ structure,
contained in \cite{P3}, Theorem 1.1, is described as follows.
Let $OH\equiv OH(I)$ be a Hilbert space equipped with the
$OH$ structure for a fixed index set $I$ and $x\in\bm_n(OH)$.
Then
$$
\|x\|_{\bm_n(OH)}=\|(x_{ij},x_{kl})\|^{1/2}_{\bm_{n^2}}
$$
where the enumeration of the entries of the numerical matrix on the l.h.s. is
as that in (\ref{dual}).

\subsection{Tensor products between operator spaces}

Let $E$, $F$ be operator spaces, and $E\otimes F$
denotes its algebraic tensor product. One can consider on $\bm_n(E\otimes F)$
several norms as above. Let $u\in\bm_n(E\otimes F)$; the norm $\|u\|_\wedge$ 
is defined as
$$
\|u\|_\wedge:=\inf\{\|\a\|\|v\|\|w\|\|\b\|\}
$$
where the infimum is taken on all the decompositions
$$
u=\sum\a_{ik}v_{ij}\otimes w_{kl}\b_{jl}
$$
with 
$\a\in\bm_{n,pq}$, $v\in\bm_p(E)$, $w\in\bm_q(F)$, $\beta\in\bm_{pq,n}$;
$\|u\|_{\vee}$ is the norm determined by the inclusion 
$E\otimes F\subset\gb(\ch\otimes\ck)$ if $E\subset\gb(\ch)$,
$F\subset\gb(\ck)$ (the last caracterization does not depend on the specific 
realization of $E$, $F$ as 
concrete operator spaces). The completions of these
tensor products are denoted respectively by 
$E\otimes_{max}F$, $E\otimes_{min}F$ and are 
referred as {\it projective} and {\it spatial} tensor product respectively,
these tensors are themselves
operator spaces (see \cite{ER2}). Also of interest is
the following complete identification $\bk_I(E)=E\otimes_{min}\bk_I$
for each operator space $E$.\\

\vskip 0.1cm

The projective tensor product allows one to describe the predual of a 
$W^*$-tensor product in terms of the preduals of its individual factors.
Namely, let $N$, $M$ be $W^*$-algebras, then the predual
$(N\overline\otimes M)_*$ is completely isomorphic to 
the projective tensor product
$N_*\otimes_{max}M_*$. The detailed proof of the above results can be 
found in \cite{ER2} Secion 3.\\

\vskip 0.1cm

The Haagerup tensor product $E\otimes_h F$ between the operator spaces
$E$, $F$, are also of interest here. It is defined as
the completion of the algebraic tensor product $E\otimes F$ under
the following norm for $u\in\bm_n(E\otimes F)$:
$$
\|u\|_h:=\inf\{\|v\|\|w\|\}
$$
where the infimum is taken on all the decompositions
$$
u_{ij}=\sum_{k=1}^p v_{ik}\otimes w_{kj}
$$
with 
$v\in\bm_{np}(E)$, $w\in\bm_{pn}(F)$, $p\in\bn$.

\subsection{The metrically nuclear maps}

The class of the metrically nuclear maps $\gd(E,F)$ between operator
spaces $E$, $F$, has been introduced and studied in \cite{ER4}. They are
defined as
$$
\gd(E,F):=E^*\otimes_{max}F/\ker\cx
$$
where $\cx$ is the map (\ref{mappa}) which is in this case
a complete quotient map. The metrically nuclear norm
is just the quotient one, see \cite{F}, Theorem 2.3. 
Another (more concrete) description
of the metrically nuclear operator has been given in \cite{F} at the
same time and independently where also a geometrical characterization
(Definition 2.6) has been made. All the spaces $\gd(E,F)$ are
themselves operator spaces which are complete if the range space
$F$ is complete. Moreover the metrically nuclear maps satisfy
the ideal property, see \cite{F} Proposition 2.4.

\section{The non-commutative vector-valued $L^p$ spaces}

Through the interpolation technique relative to the operator spaces, 
the vector-valued non-commutative 
$L^p$ spaces were introduced and intensively studied by Pisier 
\cite{P2,P4}. In this Section we summarize some of the properties
of the Pisier non-commutative $L^p$ spaces.\\ 
We consider together an operator space $E$, the
compatible couple of operator spaces
$$
\left(S_\infty(H)\otimes_{min}E,S_1(H)\otimes_{max}E\right)
$$
where $S_1(H)$, $S_\infty(H)$ are the trace class and the
class of all the compact operators acting on the Hilbert space $H$.
Then the vector-valued non-commutative $L^p$ spaces $S_p[H,E]$
can be defined as the interpolating spaces relative to the
above compatible couple:
\begin{equation}
\label{sp}
S_p[H,E]:=\left(S_\infty(H)\otimes_{min}E,S_1(H)\otimes_{max}E\right)_\th
\end{equation}
where $\th=1/p$.
If $F\subset\gb(L)$ is another operator space and $a,b\in S_{2p}(L)$,
we consider the linear map $M_{a,b}$
between $\gb(L)$, $S_{p}(L)$ defined as $M_{a,b}x:=axb$ and the map
${\widetilde M}_{a,b}$ between $S_p[H,E]\otimes_{min}\gb(L)$, 
$S_p[H,E]\otimes_{min}S_{p}(L)$ associated to $I\otimes M_{a,b}$. 
Then the following theorem allows us to compute the norms on all 
$\bm_n(S_p[H,E])$.
\begin{thm}$($\cite{P2}, Th\'eor\`eme 2$)$
\label{norma}
Let $1\leq p<+\infty$.
\itm{i} For $u\in S_p[H,E]$ one gets
$$
\|u\|_{S_p[H,E]}=\inf\{\|a\|_{S_{2p}(H)}
\|v\|_{S_\infty[H,E]}\|b\|_{S_{2p}(H)}\}
$$
where the infimum is taken on all the decomposition 
$u=(a\otimes I_E)v(b\otimes I_E)$ with $a,b\in S_{2p}(H)$ and
$v\in S_\infty[H,E]$.
\itm{ii} If $F\subset\gb(L)$ is another operator space and
$u\in S_p[H,E]\otimes_{min}F$, then its norms is given by
$$
\|u\|=\sup\{\|{\widetilde M}_{a,b}u\|_{S_p[H\otimes L,E]}\}
$$
where the supremum is taken for all $a$, $b$ in the unit ball of
$S_{2p}(L)$.
\end{thm}
As it has been described in \cite{P4} Theorem 1.1, the non-commutative
vector-valued $L^p$ spaces can be viewed as Haagerup tensor products
as follows
$$
S_p[H,E]=R(1-\th)\otimes_hE\otimes_h\overline {R(\th)}
$$
where $R(\th):=(H_r,H_c)_{\th}$ with $\th=1/p$.
In particular we have for $p=2$
$$
S_2[H,E]=OH\otimes_hE\otimes_h\overline {OH}.
$$
If the Hilbert space is kept fixed, we write $S_p[E]$
instead to $S_p[H,E]$.\\
 
We conclude with a result quite similar to that contained
in \cite{ER2} Proposition 3.1 which will be useful
in the following. Let $H$ be a Hilbert space
of (Hilbert) dimension given by the index set
$I$. Making the identification $H\equiv\ell^2(I)$ we get the following
\begin{prop}
\label{matr}
An element $u$ in $S_p[H,E]\subset\bm_I(E)$ satisfies $\|u\|_{S_p[H,E]}<1$
iff there exists elements $a,b\in S_{2p}(H)\subset\bm_I(E)$ with and
$v\in\bm_I(E)$ with $\|a\|_{S_{2p}(H)}=\|b\|_{S_{2p}(H)}=1$ and  
$\|v\|_{\bm_I(E)}<1$ such that
$$
u=avb.
$$
Furthermore one can choose $v\in\bk_I(E)$.
\end{prop}
\begin{pf}
By Theorem \ref{norma}, part $(i)$, we have to prove only 
the if part of the statement.\\
Suppose that $u\in\bm_I(E)$ can be written as $u=avb$ as above, 
and let $\eps>0$ be fixed. Then there exists $F(\eps)$
such that
\begin{align*}
\|a^{F_1}-a^{F_2}\|_{S_{2p}[E]}&\\
&\leq\eps/(3\|v\|_{\bm_I(E)})\\
\|b^{G_1}-b^{G_2}\|_{S_{2p}[E]}&
\end{align*}
whenever $F_1,F_2,G_1,G_2\supset F(\eps)$ are finite
subsets of $I$.
Consider the finite subsets $F,G,\widehat F,\widehat G\subset I$, we get
\begin{align*}
\|a^Fvb^G-a^{\widehat F}vb^{\widehat G}\|_{S_p[E]}
\leq&\|(a^F-a^{F\wedge\widehat F})vb^G\|_{S_p[E]}\\
+&\|(a^{\widehat F}-a^{F\wedge\widehat F})vb^{\widehat G}\|_{S_p[E]}\\
+&\|a^{F\wedge\widehat F}v
(b^{G\vee\widehat G}-b^{G\wedge\widehat G})\|_{S_p[E]}.
\end{align*}
Now, if $F,G,\widehat F,\widehat G\supset F(\eps)$,
we obtain 
$$
\|a^Fvb^G-a^{\widehat F}vb^{\widehat G}\|_{S_p[E]}\leq\eps.
$$
Then $\{a^Fvb^G\}$, $F,G\subset I$ finite subsets, is a Cauchy net in
$S_p[E]$ which converges to an element of $S_p[E]$ which must coincide
with $u$.
\end{pf}

\section{The (operator) $p$-factorable maps}

In this Section we introduce for $1\leq p<+\infty$, a class $\gf_p(E,F)$
of linear maps between operator spaces $E$, $F$ which are limits of finite
rank maps so $\gf_p(E,F)\subset\gk(E,F)$. These maps are obtained 
considering operators arising in a natural way from the Pisier 
non-commutative vector-valued $L^p$ spaces $S_p[H,E]$;
then we obtain the $\gf_p(E,F)$
considering those operators which factor via such maps. Although the
case with $p=+\infty$ presents no complications, to simplify,
we deal only with the cases with $1\leq p<+\infty$ (where,
as usual, the conjugate exponent of $p=1$ is $q=+\infty$).\\

We start with some elementary technical results.\\ 

\vskip 0.1cm

In the sequel we indicate with $\underline x$ any element of $\bc^n$,
$n$ arbitrary; if $\{a_i\}_{i=1}^n\subset E$, we denote the numerical
sequence $\{\|a_1\|,\dots,\|a_n\}$ simply with $\underline a$.\\

Let $\underline x:=(x_1,x_2)$, $\underline y:=(y_1,y_2)$ be four positive number
and $1\leq p\leq+\infty$ with $q={p\over p-1}$ the conjugate exponent, we have
\begin{lem}
\label{dir1}
$$
\|\underline x\|_p\|\underline y\|_q
\leq C(p)(\|\underline x\|_2^2+\|\underline y\|_2^2)
$$
where $C(p)$ is a constant which is equal to $1/2$ if $p=2$
and is greater than $1/2$ for $p\neq2$.
\end{lem}
\begin{pf}
The case with $p=2$ is elementary, so we treat only the other
cases. We get
\begin{align*}
\|\underline x\|_p\|\underline y\|_q\leq{1\over2}
(\|\underline x\|^2_p+\|&\underline y\|^2_q)\leq
{1\over2}(\|\underline x\|^2_1+\|\underline y\|^2_1)=\\
{1\over2}(x_1^2+x_2^2+2x_1x_2+y_1^2&+y_2^2+2y_1y_2)\leq
\|\underline x\|_2^2+\|\underline y\|_2^2
\end{align*}
so $C(p)\leq1$. If $1\leq p<2$ we consider the following case:
$x_1=x_2=1$, $y_1=y$, $y_2=0$. An elementary computation tell
us that $C(p)>{1\over2}$ if $p\neq2$.
\end{pf}
Now we consider the following situation where
$1\leq p<+\infty$ and $q$ is the conjugate exponent of $p$.\\

Let $A_i\in\gam(E,S_p^*),\ i=1,2$ be completely bounded maps and
consider the linear map between $E$, $S_p^*$ given by
$$
Ax:=A_1x\oplus A_2x
$$
(where we have kept fixed any identification 
$H\equiv\ell^2\cong H\oplus H$).
At the same way, let $b_i\in S_p[E],\ i=1,2$ and consider the element 
$b\in\bm_\infty(E)$ given by
$$
b:=b_1\oplus b_2.
$$
We get the following
\begin{lem}
\label{dir2}
\itm{i} $A\in\gam(E,S_p^*)$ with $\|A\|_{cb}\leq\|\underline A\|_q$,
\itm{ii} $b\in S_p[E]$ with $\|b\|_{S_p[E]}\leq\|\underline b\|_p$.
\end{lem}
\begin{pf}
The case with $p=1$ in $(i)$ is easy and is left to the reader.
For the other cases in $(i)$, taking into
account \cite{P4} Corollary 1.3 and Theorem \ref{norma}, part $(ii)$,
we compute for $n$ integer, $u,v\in(\bm_n)_1$ and $\|x\|_{\bm_n(E)}<1$
\begin{align*}
\|\widetilde M_{u,v}Ax\|^q_{S_q((H\oplus H)\otimes\bc^n)}&=
\|\widetilde M_{u,v}A_1x\|^q_{S_q(H\otimes\bc^n)}+
\|\widetilde M_{u,v}A_2x\|^q_{S_q(H\otimes\bc^n)}\\
\leq\|A_1x\|^q_{\bm_n(S_q(H))}&+\|A_2x\|^q_{\bm_n(S_q(H))}\leq
\|A_1\|^q_{cb}+\|A_2\|^q_{cb}.
\end{align*}
Taking the supremun on the left, first on the unit balls of $(\bm_n)_1$, $E$
and then on $n\in\bn$, we obtain the assertion again by Theorem \ref{norma}.\\
The proof of part $(ii)$ follows at the same way.
\end{pf}
Let any identification $\ell^2\equiv H\cong\oplus_{i=1}^NH$ be fixed and
$1\leq p<+\infty$. Suppose that we have a sequence 
$\{A_i\}_{i=1}^N\subset\gam(S_p(H),E)$. We define a linear operator
$A:S_p\left(\oplus_{i=1}^NH\right)\to E$ as follows. Let 
$x\in S_p\left(\oplus_{i=1}^NH\right)$, first we cut the off-diagonal part
of $x$, so can define
\begin{equation}
\label{off}
Ax:=\sum_{i=1}^NA_iP_ixP_i
\end{equation}
where $P_j$ is the ortogonal projection corresponding to the $j$-subspace
in the direct sum $\oplus_{i=1}^N$.
We have the following 
\begin{lem}
\label{ff}
The map $A:S_p\left(\oplus_{i=1}^NH\right)\to E$ defined as above
is completely bounded and
$$
\|A\|_{cb}\leq\|\underline A\|_q
$$
where $q$ is the conjugate exponent of $p$.
\end{lem}
\begin{pf}
It is easy to note that $A$ is bounded as an operator between
$S_p\left(\oplus_{i=1}^NH\right)$, $E$. Hence, to compute its completely
bounded norm, it is enough to pass to the transpose map
$A^*:E^*\to S_p\left(\oplus_{i=1}^NH\right)^*$ which have the same
form as that in the preceding lemma. So a similar calculation shows
that
$$
\|A\|_{cb}\equiv\|A^*\|_{cb}\leq\|\underline A^*\|_q
\equiv\|\underline A\|_q
$$
which is the assertion.
\end{pf}
Now we are ready to define the classes of factorable maps between
operator spaces.
\begin{defin}
Let $E$, $F$ be operator spaces and $1\leq p<+\infty$.\\
A linear map $T:E\to F$ will be called $p$-{\it factorable} if there
exists a Hilbert space $H$ and elements $b\in S_p[H,E^*]$, 
$A\in\gam\left(S_p(H),F\right)$ such that $T$ factorizes according to

\begin{figure}[hbt]
\begin{center}
\begin{picture}(350,120)
\put(100,20){\begin{picture}(150,100)
\thicklines
\put(33,75){\vector(1,0){84}}
\put(28,65){\vector(1,-1){33}}
\put(83,32){\vector(1,1){33}}
\thinlines
\put(20,70){$E$}
\put(58,20){$S_p[H]$}
\put(120,70){$F$}
\put(70,83){$T$}
\put(32,40){$B$}
\put(110,40){$A$}
\end{picture}}
\end{picture}
\end{center}
\end{figure}

where, in the above diagram, $B=\cx(b)\in\gam\left(E,S_p(H)\right)$, see 
\cite{P4}, Lemma 3.15.\\
We also define for a $p$-factorable map $T$
$$
\f_p(T):=\inf\{\|A\|_{cb}\|b\|_{S_p[H,E^*]}\}
$$
where the infimum is taken on all the factorization for $T$ as above.
The class of all the $p$-factorable maps between $E$, $F$ will be
denoted as $\gf_p(E,F)$.
\end{defin}
\begin{rem}
As the linear map $\cx(b)$ is norm limit of finite rank maps, hence 
has separable range, without loss of generality 
we can reduce ouselves in the above Definition to consider $H\equiv\ell^2$
and omit the dependance on $H$ in the following if it is not 
otherwise specified.
\end{rem}
\begin{rem}
\label{11}
We have $\gf_p(E,F)\subset\gk(E,F)$ as 
$\cx(\bk_I(E^*))\subset\gk(E,\bk_I)$ where, as usual, $\cx$ is the map defined
in (\ref{mappa}).
\end{rem}
Now we show that the sets of maps $\gf_p(E,F)$ are indeed 
quasi-normed linear spaces.
\begin{prop}
\label{quasi}
$\left(\gf_p(E,F),\f_p\right)$ is a quasi-normed vector space for
each $1\leq p<+\infty$. Moreover $\left(\gf_2(E,F),\f_2\right)$
is a normed vector space.
\end{prop}
\begin{pf}
If $T=A\cx(b)\in\gf_p(E,F)$, one has 
$\|T\|\leq\|A\|_{cb}\|b\|_{S_p[E^*]}$
and, taking the infimum on the right, one obtains $\|T\|\leq\f_p(T)$ so
$\f_p(T)$ is nondegenerate. We have only to prove the (generalized)
triangle inequality for $\f_p$. Let $T_i\in\gf_p(E,F)$, $i=1,2$ and 
$\eps>0$ be fixed and choose $A_i$, $b_i$ such that
$$
\|A_i\|_{cb}=\|b_i\|_{S_p[E^*]}=\sqrt{\f_p(T_i)(1+\eps)}.
$$
We consider the linear map $A:S_p\to F$ defined as in (\ref{off})
where any decomposition of $\ell^2\equiv H\cong H\oplus H$ 
has been considered. We also
consider under this above decomposition of $H$, $b:=b_1\oplus b_2$.
Applying Lemma \ref{ff}, Lemma \ref{dir2} 
we obtain $\|A\|_{cb}\leq\|\underline A\|_q$,
$\|b\|_{S_p[E^*]}=\|\underline b\|_p$; hence $T_1+T_2=A\cx(b)$.\\
By Lemma \ref{dir1}, we get
\begin{align*}
&\f_p(T_1+T_2)\leq\|A\|_{cb}\|b\|_{S_p[E^*]}\leq\\
C(p)(\|A_1&\|^2_{cb}+\|b_1\|^2_{S_p[E^*]}
+\|A_2\|^2_{cb}+\|b_2\|^2_{S_p[E^*]})\leq\\
&2C(p)(1+\eps)(\f_p(T_1)+\f_p(T_2))
\end{align*}
and the proof now follows as $\eps$ is arbitrary.
\end{pf}
Actually $\left(\gf_2(E,F),\f_2\right)$ is a quasi-normed complete
linear space
of linear maps between $E$, $F$. The case with $p=2$ gives rise to 
Banach spaces of maps.
Following the considerations contained in Section 2, one has 
for $T\in\gf_p(E,F)$ a summation 
$$
T=\sum_{i\in\bn}f_i(\cdot)y_i
$$
where $\{f_i\}\subset E^*$ $\{y_i\}\subset F$. It is easy to see that such a 
summation is unconditionally convergent in the norm topology of $\gb(E,F)$.
Moreover, if $T\in\gf_p(E,F)$, then $T$ is completely bounded, 
see Remark \ref{11}.

\section{Ideals between operator spaces}

As we have already mentioned, one can point out in a natural way 
the properties which characterize classes of operator ideals 
also in the non-commutative functional analysis that is
in operator spaces setting. Examples of such operator ideals have been
studied in \cite{ER3,ER4,ER5,F,1P,P2,P4} 
where it has been shown that some of these spaces of maps
also have a natural operator space
structure themselves. In this section we start with these definitions
and show that the factorable maps $\gf_p$ are examples of quasi-normed 
(complete) operator ideals in operator spaces setting. Moreover
if $p=2$ we obtain another example of Banach operator ideal. In order
to do this we follow the strategy of the celebrated monograph \cite{Pi1} of
Piestch.\\

We indicate with $\gam$ the classes of all the completely bounded maps, 
that is $\gam(E,F)$ is just the space of the completely bounded maps between
the operator spaces $E$, $F$.\\

The following definition is our startpoint.
\begin{defin}
A subclass $\gi\subset\gam$ will be said an {\it operator ideal}
if
\itm{i} $I_1\in\gi$ where $1$ is the $1$-dimensional space,
\itm{ii} $\gi(E,F)$ is a linear space for every operator spaces $E$, $F$,
\itm{iii} $\gam\gi\gam\subset\gi$ (ideal property).\\

Moreover if there exists quasi-norms $\f$ (\cite{Ko}, Section 15.10) such that 
\itm{a} $\f(I_1)=1$,
\itm{b} $\f(S+T)\leq \kappa(\f(S)+\f(T)),\quad \kappa\geq1$,
\itm{c} If $T\in\gam(E_0,E)$, $S\in\gi(E,F)$, $R\in\gam(F,F_0)$
then 
$$
\f(RST)\leq\|R\|_{cb}\f(S)\|T\|_{cb}
$$
with each $(\gi(E,F),\f)$ complete as topological vector space,
we call $(\gi,\f)$ a {\it quasi normed} or {\it Banach} operator ideal 
according with $\kappa>1$ or $\kappa=1$ respectively.
\end{defin}
Now we show that the $p$-factorable maps $(\gf_p,\f_p)$ are quasi-normed
operator ideals and in particular $(\gf_2,\f_2)$ is a Banach operator ideal.
For reader's convenience we split up the proof in two propositions.
\begin{prop}
Let $E_0$, $E$, $F$, $F_0$, be operator spaces and $T:E_0\to E$,
$S:E\to F$, $R:F\to F_0$ linear maps.
If $T\in\gam(E_0,E)$, $S\in\gf_p(E,F)$, $R\in\gam(F,F_0)$
then $RST\in\gf_p(E_0,F_0)$ and
$$
\f_p(RST)\leq\|R\|_{cb}\f_p(S)\|T\|_{cb}.
$$
\end{prop}
\begin{pf} If $S\in\gf_p(E,F)$ and $\eps>0$, by Proposition \ref{matr}, 
there exists $a,b\in (S_{2p})_1$, $f\in\bm_\infty(E^*)$,
$A\in\gam(S_p,F)_1$ such that 
$\|f\|_{cb}\leq\f_p(S)+{\displaystyle\eps\over{\|R\|_{cb}\|T\|_{cb}}}$ 
and $S=A\cx(af(\cdot)b)$ where $\cx$ is the map given in (\ref{mappa}).
Now 
$$
RSTx=RAa(f(Tx)b
$$ 
where $f\circ T\in\bm_\infty((E_0)^*)$ and 
$\|f\circ T\|_{cb}\leq\|f\|_{cb}\|T\|_{cb}$. Then we obtain
$RST\in\gf_p(E_0,F_0)$ and 
$$
\f_p(RST)\leq\|R\|_{cb}\|f\circ T\|_{cb}
\leq\|R\|_{cb}\|f\|_{cb}\|T\|_{cb}
\leq\|R\|_{cb}\f_p(S)\|T\|_{cb}+\eps
$$
and the proof follows.
\end{pf}
\begin{prop} 
$(\gf_p(E,F)),\f_p)$ is a complete quasi-normed space.
\end{prop}
\begin{pf} We have already proved in Proposition \ref{quasi} that
$(\gf_p(E,F)),\f_p)$ is a quasi-normed linear space for each $p$.
According with \cite{Pi1}, Section 6.2 it is enough to show that an 
absolutely summable sequence 
$\{T_i\}\subset(\gf_p(E,F),\f_p)$ is summable in $(\gf_p(E,F),\f_p)$. 
We compute $r:={1\over 2+\log_2C(p)}\leq1$ as $C(p)\geq{1\over2}$.
Let $\{T_i\}$ be an absolutely summable sequence 
(i.e. ${\displaystyle \sum_{i=1}^{+\infty}\f_p(T_i)^r<+\infty}$)
where
$T_i=A_i\cx(b_i)$ for sequences $\{A_i\}\subset\gam(S_p,F)$, 
$\{b_i\}\subset S_p[E^*]$, 
\begin{align*}
&\|A_i\|_{cb}\leq\left(\f_p(T_i)+{\eps\over2^{i+1}}\right)^{{p-1\over p}}\\
&\|b_i\|_{S_2[E^*]}\leq\left(\f_p(T_i)+{\eps\over2^{i+1}}\right)^{{1\over p}}
\end{align*}
and $H\cong\ell^2$ as usual; moreover, since $r\leq1$, 
there is a constant $K>0$ such that 
$$
\sum_{i=1}^{+\infty}\f_p(T_i)\leq K\sum_{i=1}^{+\infty}\f_p(T_i)^r<+\infty.
$$
We fix in the following any identification 
$H\equiv\ell^2\cong\oplus_{i=1}^{+\infty}H$. We define
$$
b:=\oplus_{i=1}^{+\infty}b_i\in\bm_\infty(E^*).
$$
It easy to show that 
$b\in S_p\left[\oplus_{i=1}^{+\infty}H,E^*\right]$ as
norm limit of the sequence 
$$
\s_N:=b_1\oplus\dots\oplus b_N\oplus0\oplus\dots
$$ 
thanks to
\begin{align*}
&\|\s_N\|^p_{S_p[\oplus_{i=1}^{+\infty}H,E^*]}
=\sum_{i=1}^N\|b_i\|^p_{S_p[H]}\\
\leq&\sum_{i=1}^N\f_p(T_i)+{\eps\over2}
\leq K\sum_{i=1}^{+\infty}\f_p(T_i)^r+{\eps\over2},
\end{align*}
see \cite{P4}, Corollary 1.3.\\
Moreover we can define, as in the preceding Section,
linear maps $A_N:S_p\left(\oplus_{i=1}^NH\right)\to F$
as in (\ref{off}), which are
completely bounded and satisfy, by Lemma \ref{ff}
$$
\|A_N\|_{cb}^q\leq\sum_{i=1}^{+\infty}\|A_i\|_{cb}^q
\leq\sum_{i=1}^{+\infty}\f_p(T_i)+{\eps\over2}
\leq K\sum_{i=1}^{+\infty}\f_p(T_i)^r+{\eps\over2}.
$$
By these consideration one can easily shows that the direct limit
${\displaystyle \lim_{\longrightarrow} A_N}$ defines a bounded map on
${\displaystyle\bigcup_NS_p\left(\oplus_{i=1}^NH\right)}$ which 
uniquely extends to a completely bounded map 
$A\in\gam\left(S_p\left(\oplus_{i=1}^{+\infty}H\right),F\right)$
as ${\displaystyle\bigcup_NS_p\left(\oplus_{i=1}^NH\right)}$
is dense in $S_p\left(\oplus_{i=1}^{+\infty}H\right)$.
So we have 
\begin{align*}
&b\in S_p\left[\oplus_{i=1}^{+\infty}H,E^*\right]\cong
S_p\left[H,E^*\right],\\
&A\in\gam(S_p\left(\oplus_{i=1}^{+\infty}H\right),F)\cong
\gam(S_p(H),F).
\end{align*}
Finally, if one defines $T:=A\cx(b)$, then $T\in\gf_p(E,F)$
and
$$
\f_p(T-T_N)\leq
\left(\sum_{i=N+1}^{+\infty}\f_p(T_i)^r+{\eps\over2}\right)^{{p-1\over p}}
\left(\sum_{i=N+1}^{+\infty}\f_p(T_i)^r+{\eps\over2}\right)^{{1\over p}}
<\eps
$$
if $N$ is big enough, that is $T$ is summable in $\gf_p(E,F)$. 
Moreover, by \cite{Pi1}, 6.1.9, 6.2.4,
one also gets
\begin{align*}
\f_p(T)^r&\leq\liminf_N\f_p(T_N)^r\\
&\leq2\liminf_N\sum_{i=1}^N\f_p(T_N)^r\\
&=2\sum_{i=1}^{+\infty}\f_p(T_N)^r<+\infty.
\end{align*}
\end{pf}
Summarizing we have the following
\begin{thm}
$(\gf_p,\f_p)$, $1\leq p<+\infty$ are quasi-normed operator ideals.
Moreover $(\gf_2,\f_2)$ is a Banach operator ideal.
\end{thm}
\begin{pf}
The proof of the first part immediately follows collecting the 
results contained in the preceding Section and in the last two Propositions.
The case with $p=2$ follows by observing that $C(2)=1/2$ so $\kappa=1$ that is
$\f_2$ is a norm and $(\gf_2(E,F),\f_2)$ is complete.
\end{pf}
\begin{rem}
On each of $\gf_p(E,F)$ there is a $r$-norm $A_r$ equivalent to $\f_p$ given by
$$
A_r(S):=\inf\{\left(\sum_{i=1}^n\f_p(S_i)^r\right)^{{1\over r}}\}
$$
where the infimum is taken on all the decomposition $S=\sum_{i=1}^n S_i$
with $S_i\in\gf_p(E,F)$, see \cite{Pi1}, 6.2.5.
\end{rem}

\section{A geometrical description}

Analogously to the metrically nuclear operators, see \cite{F},
we give a suitable geometrical description for the range of a $p$-factorable
injective map.\\

\vskip 0.1cm

We start with an absolutely convex set $Q$ in an operator space
$E$ and we indicate with $V$ its algebraic span. Consider a sequence
$\cq\equiv \{Q_n\}$ of sets such that

\begin{itemize}
\item[(i)] $Q_1\equiv Q$ and every $Q_n$ is an absolutely convex absorbing
set of $\bm_n(V)$ with $Q_n\subset\bm_n(Q)$;
\item[(ii)] $Q_{m+n}\cap(\bm_m(V)\oplus\bm_n(V))=Q_m\oplus Q_n$;
\item[(iii)] for $x\in Q_n$ then 
$x\in\lambda Q_n$ implies $bx\in\lambda Q_n, xb\in\lambda Q_n$
where $b\in(\bm_n)_1$.\\
\end{itemize}

We say that a (possibly) infinite matrix $f$ with entries in the algebraic dual
$V'$ of $V$ has {\it finite $\cq$-norm} if
$$
\|f\|_\cq\equiv\sup\{\|f^\Delta(q)
\|:q\in Q_n;\ n\in\bn;\ \Delta\}<+\infty
$$
where $f^\Delta$ indicates an arbitrary 
finite truncation corresponding to the finite set $\Delta$; 
the enumeration of the entries of the numerical matrix 
$f^{\Delta}(q)$ is as that in (\ref{dual}).

\begin{defin}
\label{fatt}
An absolutely convex set $Q\subset E$ is said to be
{\it $(p,\cq)$-factorable}, $1\leq p<+\infty$ (where $\cq$ is a fixed 
sequence as above)
if there exists matrices $\a,\b\in S_{2p}$ and a (possible infinite)
matrix $f$ of linear functionals as above with $\|f\|_\cq<+\infty$
such that, if $x\in Q_n$, one has
\begin{equation}
\label{ultra}
\|x\|_{\bm_n(E)}\leq C\|\a f(x)\b\|_{\bm_n(S_p)}
\end{equation}
In the case corresponding to $p=2$ we call a $(2,\cq)$-factorable
set simply $\cq$-factorable and omit the dependence on $\cq$ if it causes
no confusion.
\end{defin}

One can easily see that a $(p,\cq)$-factorable 
set is relatively compact,
hence bounded in the norm topology of $E$ and therefore $V$, together the
Minkowski norms determinated by the $Q_n$'s on $\bm_n(V)$, is a (not necessarily
complete) operator space. As in the metrically nuclear case described in
\cite{F}, Section 2, and the case of completely summing maps described in
\cite{P4} Remark 3.7, one can reinterpret the above definition as
a factorization condition.
\begin{prop}
Let $E\subset\gb(\ch)$ be a (concrete) operator space, $Q\subset E$
a $(p,\cq)$-factorable absolutely convex set and $V$
its algebraic span.
Then the canonical immersion $V\stackrel{i}{\hookrightarrow}E$ is a 
$p$-factorable map when $V$ is equipped with the operator space 
structure determined by the sequence $\cq$
\itm{i} if there exists a completely bounded projection $P:\gb(\ch)\to E$
when $p\neq2$;
\itm{ii} without any other condition on $E$ if $p=2$.
\end{prop}
\begin{pf}
According to the above Definition, there exists 
matrices $\a,\b\in S_{2p}$ and a 
matrix $f$ of linear functionals with $\|f\|_\cq<+\infty$
satisfying the property described above.
We define $b:=\a f\b$ so $b\in S_2[V^*]$.
Moreover we can consider $W:=\overline {\a f(V)\b}^{S_p}$ and define on  
$W$ a linear map $A:W\to E$ in the following way $Ax:=v$ if $x=\a f(v)\b$.
This map extends firstly to all of $W$ and successively to a completely 
bounded map between $S_p$ and $\gb(\ch)$ by 
the celebrated Arveson-Wittstock-Hahn-Banach Theorem, \cite{W1} (see also
\cite{W3,Pa2}). Then we obtain $i=PA\cx(b)$
which is $p$-factorable as $PA$ is a completely bounded map between $S_2$ and
$E$, see (\ref{ultra}). In the case with $p=2$ we can extend $A$ to
all of $S_2$ if one define $Ax:=0$ on $W^\perp$. 
Being $S_2\cong OH$ a
homogeneous Hilbertian operator space, see \cite{P3}, Proposition 1.5, 
$A$ is completely bounded and the proof is now complete.
\end{pf}
We now consider an injective completely bounded operator
$T:E\to F$ and the sequence $\cq_T$ given by
$$
\cq_T=\{T_n(\bm_n(E)_1)\}_{n\in\bn}.
$$
for such sequences the properties (i)--(iii) in the beginning 
are automatically satisfied and if
$T(E_1)$ is $(p,\cq_T)$-factorable we call it simply $(p,T)$-factorable
and indicate the $\cq_T$-norm of a matrix of functional 
$f$ by $\|f\|_T$.\\

As it happens in some interesting well-known cases
(compact operators, nuclear and metrically nuclear maps), also for
the class of $p$-factorable maps we have a description in
terms of geometrical properties (i.e. the shape) of the range of such maps.
\begin{prop} 
\label{shape}
Let $E$, $F$ be operator spaces with $F\subset\gb(\ch)$ and $T:E\to F$ a
completely bounded injective operator. Then $T\in\gf_p(E,F)$ iff $T(E_1)$ is a
$(p,T)$-factorable set in $F$
\itm{i} if there exists a completely bounded projection $P:\gb(\ch)\to F$
when $p\neq2$;
\itm{ii} without any other condition on $F$ if $p=2$.
\end{prop}
\begin{pf}
It is easy to verify that, if $T\in\gf_2(E,F)$ is injective, one can
write for $T$ a decomposition
$$
T=A\cx(\a f(\cdot)\b).
$$
Then $\a$, $f\circ T^{-1}$, $\b$, allow us to say that
$T(E_1)$ is $(p,T)$-factorable with $\|A\|_{cb}\leq C$, $C$
is the constant appearing in (\ref{ultra}). 
Conversely, if $T(E_1)$ is a $(p,T)$-factorable set 
in $F$, with
$\a$, $\b$, $f$ as in the Definition \ref{fatt}, then we can consider
the matrix $f\circ T\in\bm_\infty(E^*)$ and define $A:=S_p\to F$ as in the
proof of the above proposition; so we obtain for $T$ the factorization
$T=A\cx(b)$ where $b:=\a(f\circ T)\b\in S_p[E^*]$.
\end{pf}
As in the case relative to metrically nuclear maps (\cite{F}),
the definition of a factorable set may appear rather involved; this is due
to the fact that the inclusion 
$\bm_n(V)_1\subset\bm_n(V_1)$ is strict in general;
but, for an injective completely bounded 
operator $T$ as above, the $T$-factorable
set $T(V_1)$ is intrinsecally defined in terms of $T$.
\begin{rem}
An interesting example in the cases with $p\neq2$ 
for which the above Propositions is applicable is when the image space  
is an injective $C^*$-algebra. 
\end{rem}

\vskip 0.1cm

A characterization of the completely $p$-summing maps in terms of
a factorization condition is given in \cite{P4}, Remark 3.7 and involve
ultrapowers (\cite{He}). As ultrapowers of $W^*$-algebras of the same type
seems to produce $W^*$-algebras of other type in general
(see \cite{Co}, Section II for same kind of similar questions, or \cite{Pi1},
19.3.4 for the commutative case) one can argue that in
general $\gf_p(E,F)\subsetneqq\Pi_p(E,F)$, the completely
summing maps considered by Pisier in \cite{P4}. 
It would be of interest to understand if (and when)
the above inclusion is in fact an equality. Other cases which could
involve ultrapowers are the definition of factorable maps in operator
space setting in a similar manner as that considered in \cite{Kw}.
Such a kind of factorable maps might be the true quantized counterpart 
of the factorable maps of Banach spaces case, see \cite{Kw,Pi,Pi1}
(see also \cite{ER5} for some related questions relative to the quantized case).
Following \cite{Pi1,Kw}, one could consider those maps 
which factors according to the following commutative diagram

\begin{figure}[hbt]
\begin{center}
\begin{picture}(350,120)
\put(100,20){\begin{picture}(150,100)
\thicklines
\put(33,75){\vector(1,0){94}}
\put(28,67){\vector(1,-1){36}}
\put(90,32){\vector(1,1){36}}
\thinlines
\put(20,70){$E$}
\put(61,20){$L^p(M)$}
\put(130,70){$F^{**}$}
\put(65,83){$i\circ T$}
\put(32,40){$B$}
\put(115,40){$A$}
\end{picture}}
\end{picture}
\end{center}
\end{figure}

Here $i$ is the canonical completely isometric immersion of $F$
in $F^{**}$ (\cite{B}) and $M$ is a $W^*$-algebra whose non-commutative 
measure spaces (\cite{K1}) should be equipped with suitable operator
space structures which are not yet understood in the full generality,
see \cite{P4}. The maps $A$, $B$ in the above commutative diagram
should be completely bounded.\\
A complete analysis of the above framework and specially 
the study of quantized counterparts of the results of Banach space
theory will be desirable. We hope to return about some of these
problems somewhere else.

\section{The split property for inclusions of $W^*$-algebras. An application to
Quantum Field Theory}

In this section we outline a characterization of the split property for 
inclusions of $W^*$-algebras by the $2$-factorable maps. 
The complete exposition of this result is contained in the 
separated paper \cite{F1}  
to which we remand the reader for further details.\\

We suppose that all the $W^*$-algebras considerated here 
have separable predual.
For the standard results about the theory of $W^*$-algebras see e.g.
\cite{SZ,T}.\\

\vskip 0.1cm

An inclusion $N\subset M$ of $W^*$-algebras is said to be {\it split} if there 
exists a type $I$ interpolating $W^*$-factor $F$ that is $N\subset F\subset M$.
The split property has been intensively studied (\cite{B1,D,DL2}) 
in the last years for the natural applications to Quantum Field Theory 
(\cite{B2,B3,B5}).
In \cite{B1}, canonical non-commutative embeddings $\F_i:M\to L^i(M)$,
$i=1,2$ are considered; they are constructed via a standard vector 
$\Om\in L^2(M)$ for $M$ in the following way
\begin{align}
\begin{split}
\label{tomita}
&\F_1:a\in M\to(\cdot\Om,Ja\Om)\in L^1(M)\\
&\F_2:a\in M\to\D^{1/4}a\Om\in L^2(M).
\end{split}
\end{align}
The $W^*$-algebra $M$ is supposed to act standarly on the Hilbert
space $L^2(M)$ and, in the above formulas, 
$\D$, $J$ are the Tomita's operators relative to the standard vector
$\Om$. In this way the split property is analyzed considering 
the nuclear properties of the restrictions of $\F_i$, $i=1,2$ to the
subalgebra $N$. The nuclear property and its connection with the
split property has an interesting physical meaning in Quantum Field Theory, 
see \cite{B2,B3,B5}. Following this approach, in \cite{F} the split 
property has been exactly characterized in terms of
the $L^1$ embedding $\F_1$ constructed by a fixed standard vector for $M$.
The characterization of the split property in terms of the other $L^2$
embedding $\F_2$ is established in \cite{F1}.
In this section we only schetch the proof of this characterization. 
which we report together with the results relative
to the $L^1$ embedding $\F_1$ for the sake of completness.
\begin{thm}
\label{mega}
Let $N\subset M$ be an inclusion of $W^*$-factors with 
separable preduals and $\omega\in M_*$ a faithful state. Let
$\F_i:M\to L^i(M)$, $i=1,2$ be the embeddings associated to the
state $\omega$ and given in (\ref{tomita}).\\ 
The following statements are equivalent.
\itm{i} $N\subset M$ is a split inclusion.
\itm{ii} $\F_{1\lceil N}\in\gd(N,(L^1(M))$.
\itm{ii'} The set $\{(\cdot\Omega,Ja\Omega):a\in A,\ \|a\|<1\}$ is
$\F_1$-decomposable (see \cite{F}, Definition 2.6 for this definition).
\itm{iii} $\F_{2\lceil N}\in\gf_2(N,(L^2(M))$.
\itm{iii'} The set $\{\D^{1/4}a\Om:a\in A,\ \|a\|<1\}$ is
$\F_2$-factorable.
\end{thm}
In the above theorem we are supposing that
$L^i(M)$, $i=1,2$ are endowed with the operator space structures:

\vskip 0.1cm

\begin{itemize}
\item[(a)] as the predual of $M^\circ$, the opposite algebra of $M$,
for the former,
\item[(b)] the Pisier $OH$ structure for the latter.
\end{itemize} 

\vskip 0.1cm

\begin{pf}
Some of the above equivalences are immediate (see Proposition \ref{shape})
or are contained in \cite{F} so we only deal with the remaining ones.\\
$(i)\Rightarrow(iii)$ If there exists a type $I$ interpolating factor $F$
then $\F_2$ factors according to

\newpage

\begin{figure}[hbt]
\begin{center}
\begin{picture}(350,120)
\put(100,20){\begin{picture}(150,100)
\thicklines
\put(33,75){\vector(1,0){84}}
\put(28,67){\vector(1,-1){36}}
\put(81,32){\vector(1,1){36}}
\thinlines
\put(20,70){$N$}
\put(60,20){$L^2(F)$}
\put(120,70){$L^2(M)$}
\put(70,83){$\F_2$}
\put(32,40){$\Psi_2$}
\put(115,40){$\Psi_1$}
\end{picture}}
\end{picture}
\end{center}
\end{figure}

where $\Psi_2$ arises from $S_2[N_*]$ and $\Psi_1$ is bounded, see \cite{B1}. 
Moreover $\Psi_1$ is automatically completely bounded, see \cite{P3}, 
Proposition 1.5.\\
$(iii)\Rightarrow(i)$ It is enough to show that, if $N\subset F$ with
$F$ a type $I$ factor with separable predual, 
then $\F_1$ extends to a completely positive
map $\widetilde{\F}_1:F\to L^1(M^\circ)$, see \cite{B1}, Proposition 1.1.\\
Suppose that $\F_2=A\cx(b)$ with $A$ completely bounded and $b\in S_2[N^*]$.
As $\F_2$ is normal, we can choose $b\in S_2[N_*]$ 
(\cite{F1}, Lemma 1). Hence, as $b=\a f\b$
with $f\in\bm_\infty(N_*)$, we have via (\ref{mappa}) a 
completely bounded normal
map $\r:N\to\bm_\infty\equiv F_\infty$ that is the (unique)
type $I_\infty$ factor with separable predual. So, if $N\subset F$, 
then $\r$ extends to a  
completely bounded normal map $\tilde\r:F\to\bm_\infty$ 
(\cite{F1}, Proposition 7). Computing
$\tilde f:=\cx^{-1}(\tilde\r)$ we get a completely bounded normal
extension $\widetilde{\F}:F\to L^2(M^\circ)$. Unfortunately $\widetilde{\F}$
might be not positive. We consider the binormal bilinear form 
$\f:F\otimes M^\circ\to\bc$ given by
$$
\f(x\otimes y):=(\widetilde{\F}(x),\D^{1/4}y^*\Om).
$$
This form uniquely define a bounded binormal form on all the $C^*$-algebra
$F\otimes_{max}M^\circ$ 
(\cite{F1}, Proposition 8) which can be decomposed in four binormal
positive functionals (\cite{T}). 
We take the positive part $\f_+$ and note that
$\f_+$, when restricted to $N\otimes M^\circ$  dominates $\om$    
given by
$$
\om(x\otimes y):=(x\Om,Jy\Om).
$$
If we consider the $GNS$ construction for $\f_+$, we obtain two
normal commuting representations $\pi$, $\pi^\circ$ of $F$, $M^\circ$
respectively, on a separable Hilbert space $\ch$ and a vector
$\xi\in\ch$ (cyclic for $\pi(F)\vee\pi^\circ(M^\circ)$) such that
$$
\f_+(x\otimes y)=(\pi(x)\pi^\circ(y)\xi,\xi)
$$
that is actually a $F-M$ correspondence (\cite{CJ}).
Then, as the restriction of $\f_+$ dominates $\om$, 
there exists a positive element $T\in\pi(N)'\wedge\pi^\circ(M^\circ)'$
such that, for $x\in N$, $y\in M^\circ$ one gets
$$
\om(x\otimes y)=(\pi(x)\pi^\circ(y)T\xi,T\xi).
$$
Now we define $\widetilde{\F}_1:F\to L^1(M^\circ)$ given by
$$
\widetilde{\F}_1(f):=(T\pi(f)T\pi^\circ(\cdot)\xi,\xi)
$$
which is a completely positive normal map which extends
$\F_1$.
\end{pf}
It is still unclear to the author if the split property can
be characterized via the $2$-factorable maps also for the
general case of inclusions of $W^*$-algebras with non-trivial centers. 
However in some
interesting cases such as those arising from Quantum Field Theory
we turn out to have the same characterization.
We suppose that the net $\co\to\ga(\co)$ of von Neumann algebras of local
observables of a quantum theory acts on the Hilbert space $\ch$ 
of the vacuum representation and satisfy all the usual assumptions (a priori
without the split property) which are typical in Quantum Field Theory;
$\Om\in\ch$ will be the vacuum vector which is cyclic and separating
for the net $\{\ga(\co)\}$, see e.g. \cite{Su}.
\begin{thm}
Let $\co\subset int(\widehat\co)$ be double cones in the physical
space-time and $\ga(\co)\subset\ga(\widehat\co)$ the corresponding
inclusion of von Neumann algebras of local observables.
The following assertions are equivalent.
\itm{i} $\ga(\co)\subset\ga(\widehat\co)$ is a split inclusion.
\itm{ii} The set $\{(\cdot a\Om,\Om):a\in\ga(\co)_1\}
\subset(\ga(\widehat\co)')_*$ is
$\F_1$-decomposable.
\itm{iii} The set $\{\D^{1/4}a\Om:a\in\ga(\co)_1\}
\subset\ch$ is $\F_2$-factorable.
\end{thm}
\begin{pf}
If $\ga(\co)\subset\ga(\widehat\co)$ is a split inclusion then $(ii)$
and $(iii)$ are true, see \cite{B1,F,F1}; conversely, if $(ii)$
or $(iii)$ are satisfied then $\F_1$ is extendible, see \cite{F} or
the $(iii)\Rightarrow(i)$ part of Theorem \ref{mega} 
which is also available in this case (\cite{F1}, Proposition 7).
Hence the map 
$$
\eta:a\otimes b\in\ga(\co)\otimes\ga(\widehat\co)'\to 
ab\in\ga(\co)\vee\ga(\widehat\co)'
$$
extends to a normal homomorphism of 
$\ga(\co)\overline{\otimes}\ga(\widehat\co)$ onto 
$\ga(\co)\vee\ga(\widehat\co)$. But this homomorphism 
is in fact an isomorphism by an argument exposed in \cite{B0}, pagg. 129-130.
Moreover, as $\ga(\co)\wedge\ga(\widehat\co)$ is properly infinite
(\cite{Ka}), the assertion now follows by Corollary 1 of \cite{D}.
\end{pf}


\vspace{1.5cm}

\end{document}